\documentclass[11pt,a4paper]{article}
\usepackage[utf8]{inputenc}
\usepackage{amsmath}
\usepackage{amsfonts}
\usepackage{amssymb}
\usepackage{graphicx}
\author{G.~K. Guskov, F.~I. Solov'eva}
\title{Properties of perfect transitive binary codes of length 15 and extended perfect transitive binary codes of length 16}

\begin{document}

\maketitle

\begin{quote}
{\small \noindent{\sc Abstract.} Some properties of perfect transitive binary codes of length 15 and extended perfect transitive binary codes of length 16 are presented for reference purposes.
\medskip

}
\end{quote}

The attached files contain some tab-delimited properties of perfect binary codes of length 15 and extended perfect codes of length 16. Classification of such codes can be found in \cite{OstPot}, but unfortunately the list of such codes' properties was not attached to the paper. We acknowledge that such properties had been computated by the authors in \cite{OstPot} and our results are presented for reference purposes only. The attached files are

\medskip

\textit{perfect15.txt} contains a list of properties of all perfect codes of length 15. The properties are: index of the code in classification \cite{OstPot}, rank, dimension of the kernel.

\textit{perfect16.txt} contains a list of properties of all extended perfect codes of length 16. The properties are: index of the code in classification \cite{OstPot}, rank, dimension of the kernel.

\textit{transitive15.txt} contains a list of properties of all perfect transitive codes of length 15. The properties are: index of the code in classification \cite{OstPot}, rank, dimension of the code's kernel, order of the code's automorphism group, number of codewords of weight 3 in the set $C + C$, order of the code's symmetry group.

\textit{transitive16.txt} contains a list of properties of all extended perfect transitive codes of length 16. The properties are: index of the code in classification \cite{OstPot}, rank, dimension of the code's kernel, order of the code's automorphism group, number of codewords of weight 4 in the set $C + C$, order of the code's symmetry group.

All computations have been carried out by using \textit{Magma} \cite{MagmaCitation} system.

\end{document}